\documentclass[12pt]{amsart}
\usepackage[T1]{fontenc}

\title{Un lemme combinatoire de H. B. Neumann} 

\author{Labib Haddad}
\address{120 rue de Charonne, 75011 Paris, France}
\email{labib.haddad@wanadoo.fr}

\usepackage{amssymb}
\usepackage{amsmath}
\usepackage{endnotes}

\newcommand{\su}{\subsection*}
\newcommand{\head}{\section*}
\newcommand{\noi}{\noindent}
\newcommand{\se}{\noi{\bf En effet}}
\newcommand{\cqfd}{\hfill{\bf cqfd}}

\newcommand{\Ž}{\'e}

\newcommand{\ˆ}{\`a}
\newcommand{\}{\`u}

\newcommand{\N}{\mathbb N}

\newcommand{\cal}{\mathcal}

\newcommand{\leqs}{\leqslant}

\newcommand{\geqs}{\geqslant}

\newcommand {\et}{\ \text{et}\ }

\newcommand {\pourtout}{\ \text{pour tout}\ }

\newcommand {\pourtous}{\ \text{pour tous}\ }

\newcommand{\inc}{\subset}

\newcommand{\lopar}{\noi \{$\looparrowright$ \ }

\begin{document}
\maketitle

\thispagestyle{empty}

\markboth{Labib Haddad}{Un lemme combinatoire de H. B. Neumann}

\hfill{\it Ajoutez quelquefois, et souvent effacez.}

\hfill{Boileau, L'Art po\Žtique}

\

\su{\bf R\Žsum\Ž}{\small On donne une d\Žmonstration notablement plus simple, et plus courte, du  r\Žsultat  de  H. B. Neumann qui  s'\Žnonce, sommairement, comme ceci. Pour toute partie, $A$, bien ordonn\Že d'un   semigroupe totalement ordonn\Ž, l'ensemble des produits d'un nombre fini quelconque d'\Žl\Žments de $A$ est lui-m\me bien ordonn\Ž. De plus, pour chaque $t$, il n'y qu'un nombre fini de ces produits \Žgaux \ˆ $t$.}

\

\

 On trouve dans  {\sc Neumann}, [ 3, p.209 ], un r\Žsultat  que {\sc Gonshor} a repris, dans [ 1, p.52 ], sous le titre {\sl Combinatorial lemma  on semigroups},  et qu'il \Žnonce comme suit. 
  
 \

\noi Soit $A$ une partie bien ordonn\Že de l'ensemble des \Žl\Žments strictement positifs d'un semigroupe  totalement ordonn\Ž. Alors, l'ensemble de toutes les  sommes, $a_1 + a_2 + \dots + a_q,$ d'un nombre fini quelconque d'\Žl\Žments de $A$ est bien ordonn\Ž et, de plus, pour chaque $t$ donn\Ž, il n'y a qu'un nombre fini de ces sommes qui soient \Žgales \ˆ  $t$ .

\

\noi La d\Žmonstration  de {\sc Gonshor} n'est  pas facile \ˆ suivre, pas plus que celle de {\sc Neumann}.  Elles foisonnent toutes les deux d'indices, de sous-indices et d'accents divers ! De son c\™t\Ž, {\sc Gonshor} sous-entend, sans le dire explicitement dans son \Žnonc\Ž,  que le semigroupe est commutatif ! En effet, l'op\Žration y est not\Že, $ + $, et il utilise la commutativit\Ž  dans le cas 4 de sa d\Žmonstration.

\

\noi Tandis  que {\sc Neumann} \Žtablit le r\Žsultat pour un semigroupe totalement ordonn\Ž quelconque, {\bf pas n\Žcessairement commutatif }!  Son op\Žration est une multiplication. Sa d\Žmonstration occupe plusieurs pages, (p.204 -209).
\ 

\lopar En lisant l'article [ 2 ] de Harkelroad et Gonshor, on trouvera une probable confirmation du fait que {\sc Gonshor}  entendait utiliser un semigroupe commutatif, alias ab\Žlien.\}

\

\noi Ce qui complique un peu les choses c'est que {\sc Neumann} et {\sc Gonshor} n'entendent pas tout \ˆ fait la m\me choses en parlant de semigroupe totalement ordonn\Ž.

\

\noi On va s'en tenir au choix de {\sc Neumann}, le {\bf plus g\Žn\Žral},  et donner une d\Žmonstration plus simple de son  lemme. Mais d'abord ces quelques pr\Žliminaires afin de fixer  notations et vocabulaire et d'\Žclaircir le sens des \Žnonc\Žs !

\

\head{Pr\Žparation}

\

\head{Qu'est-ce qu'un semigroupe totalement ordonn\Ž selon Neumann}

\

\noi Par d\Žfinition, un semigroupe totalement ordonn\Ž est muni d'une op\Žration binaire,  le produit, et  d'une relation d'ordre total, $\leqs$. 

\noi Le produit est associatif, autrement dit, on a 

$a(bc)= (ab)c$.

\noi Le produit et l'ordre total sont compatibles. Autrement dit,  

\noi lorsque l'on a $a <b$,  on aura $ac < bc \et ca < cb, \pourtout c$.

\noi Le plus petit de tous les  \Žl\Žments, c'est un \Žl\Žment unit\Ž, $1$.

\

On en d\Žduit, sans grands d\Žtours, les propri\Žt\Žs suivantes de simplifications \ˆ droite  et \ˆ gauche.

\noi Si $ac = bc$ alors $a = b$. 

\noi De m\me, si $ca = cb$ alors $a = b$.

\noi Donc $a = b \iff ac = bc \iff ca = cb$.

\noi Ainsi,  $ac < bc  \iff ca < cb  \iff a < b, \pourtous a, b, c$.

\noi Enfin,  $a \leqs b \iff ac \leqs bc \iff ca \leqs cb$.

\noi En particulier, $a > 1 \iff aa > a$,
 
 \noi et  $aa= 1 \iff a=1$.
 
\

\lopar En ayant adjoint un \Žl\Žment unit\Ž, $1$, plus petit que tous les autres, on simplifie consid\Žrablement les consid\Žrations de {\sc Neumann}.\}

\

\su{Classes archim\Ždiennes} On introduit la relation binaire suivante
entre les \Žl\Žments du semigroupe :
$$a\sim b \iff \text{il existe}\ m , n\in \N, \text{tels que} \ b \leqs a^m \et a\leqs b^n.$$
On v\Žrifie qu'il s'agit d'une relation d'\Žquivalence dont les classes sont des intervalles qui forment une partition du semigroupe. Pour chaque \Žl\Žment $a$ du semigroupre, on d\Žsigne par $cl(a)$ la classe d'\Žquivalence (archim\Ždienne) \ˆ laquelle il appartient, autrement dit, 

$cl(a) = \{b : b\sim a\}$,

$cl(1) = \{1\}$, en particulier.

\noi L'ensemble des classes archim\Ždiennes est totalement ordonn\Ž par la relation suivante
$$cl(a) < cl(b) \iff a^n < b, \pourtout n\in \N.$$

\su{Remarque} La classe d'un produit, $p = ab\dots c\dots d$, est \Žgale \ˆ la  classe du plus grand facteur ! En effet, supposons que tous les facteurs sont $\leqs c$. Alors $c \leqs p$ car tous les facteurs sont $\geqs 1$, d'une part. D'autre part $p \leqs c^n$ o\ $n$ est le nombre de facteurs.

\

\head{Ensembles bien ordonn\Žs. Rappels}

\

Par d\Žfinition, un ensemble {\bf  bien ordonn\Ž} est un ensemble ordonn\Ž dans lequel toute partie non vide poss\de un \Žl\Žment plus petit que tous les autres ! Un ensemble ordonn\Ž, $E$, est bien ordonn\Ž si et seulement s'il n'existe aucune suite strictement d\Žcroissante dans $E$. C'est une caract\Žrisation bien connue des ensemles bien ordonn\Žs. En voici une autre, peut-\tre moins souvent cit\Že. 

\

Soit $E$ un ensemble ordonn\Ž. Les deux \Žnonc\Žs suivants sont \Žquivalents. 

1. L'ensemble $E$ est bien ordonn\Ž.

2. De toute suite, $x_1, x_2, \dots$, dans $E$ on peut extraire une sous-suite croissante, $x_{\sigma(1)} \leqs x_{\sigma(2)}, \dots$.

\

\se, si $E$ est bien ordonn\Ž, de la suite $x_1,  x_2, \dots$, on extrait $x_{\sigma(1)}$, le plus petit terme de la suite, ayant le plus petit indice,  puis $x_{\sigma(2)}$, le plus petit petit terme de ce qui reste de la suite, ayant le plus petit indice, et on recommence ! On obtient \Žvidemment une sous-suite extraite croissante !

\

\noi R\Žciproquement, si l'\Žnonc\Ž 2 est satsifait, il n'existe \Žvidemment pas de suite strictement d\Žcroissante dans $E$. \cqfd

\

\head{Les s\Žquences}

\

\noi Soit $A$ une partie bien ordonn\Že d'un  semigroupe totalement ordonn\Ž telle que $1\notin A$. On appellera  {\bf s\Žquence} tout $q$-uplet, $s = (a, b, \dots, c)$,  o\ $a, b, \dots c$, est une suite finie  de $q$ \Žl\Žments de $A$. On dit que la s\Žquence est de longueur $q$, et on pose $|(a, b, \dots, c)| = ab\dots c$, le produit.

\

\noi Attention, le semigroupe n'\Žtant pas n\Žcessairement commutatif, on ne pourra  pas  interchanger les termes d'un produit !  Bien entendu, deux s\Žquences distinctes, de longueurs \Žgales ou  diff\Žrentes, peuvent avoir  le m\me produit !

\

\noi Il faut bien distinguer une s\Žquence, $s$, de son produit, $|s|$.  On d\Žsigne par $B$ l'ensemble de toutes les s\Žquences et par $C = \{ |s| : s \in B\}$ l'ensemble de tous les produits, de sorte que l'on a
 $$A\inc B \ , \ A\inc C \ , \ \text{et une application }\ , \ B\to C \ , \ s\mapsto |s|.$$
On dit que la s\Žquence $s$ est {\bf un repr\Žsentant} du produit $|s|$.

\

\noi On va \Žtablir le r\Žsultat suivant, par l'absurde. La premi\re partie du lemme combinatoire de {\sc Neumann} en d\Žcoulera. On essaiera d'\Žviter l'usage excessif d'indices  !

\su{Lemme A}{\sl Il n'existe aucune suite de s\Žquences,} 
$$s_1, s_2, \dots, s_n, \dots, \ \text{telles que}$$
$$|s_1| >  |s_2| > \dots > |s_n| > \cdots$$
\su{D\Žmonstration} On suppose le contraire. On suppose qu'il existe de telles suites. Soit $\cal  S$ l'ensemble de toutes ces suites. On en prend une. Pour $s_n = (a_n, b_n, \dots, c_n, \dots, d_n)$, on pose 
$$m_n = \max\{a_n, b_n, \dots, c_n, \dots, d_n\}.$$
On aura ainsi 
$$|s_n| = a_nb_n \dots c_n \dots d_n \geqs m_n  \et cl(|s_n|) = cl(m_n),$$ d'apr\s la remarque faite ci-dessus. Les $m_n$ sont des \Žl\Žments de l'ensemble bien ordonn\Ž, $A$. Soit $m$ le plus petit de ces $m_n$. La classe $cl(m)$ est la plus petite des classes  $cl(|s_n|)$. On a $cl(|s_n|) \geqs cl(|s_{n+1}|)$ puisque $|s_n| > |s_{n+1}|$. On aura donc, \ˆ partir d'un certain rang, 
$$cl(m) = cl(|s_n|) = cl(|s_{n+1}|) =cl(|s_{n+2}|) = \cdots,$$
autrement dit, tous les produits $|s_n|$ seront dans la m\me classe archim\Ždienne, $cl(m)$. En supprimant quelques uns des premiers termes de la suite, on se ram\ne au cas o\ l'on a
$$cl(m) = cl(|s_1|) = cl(|s_{2}|) =cl(|s_{3}|) = \cdots.$$
Chacun des produits $|s_n|$ est de la forme :

$|s_n] = p_nm_nr_n$ o\ $p_n$ et $r_n$ sont eux-m\mes des produits ou bien, exceptionnellement, $p_n  = 1$ ou $r_n = 1$. On a $m_n \geqs m$ pour tout $n$. On peut donc extraire une sous-suite de $p_n$ strictement d\Žcroissante, ou bien une sous-suite de $r_n$ strictement d\Žcroissante car, sinon, ou aurait $p_{n+1} \geqs p_n$ et $r_{n+1} \geqs r_n$, \ˆ partir d'un certain rang, donc
$|s_{n+1}| \geqs |s_n|$.

\

\noi Parmi toutes les suites appartenant \ˆ l'ensemble $\cal S$, on en prend une pour laquelle la classe $cl(m)$ est la plus petite possible ! Puis on d\Žsigne par $a$ le plus petit des \Žl\Žments de $A$ qui sont dans $cl(m)$. Il vient
$$|s_n| \geqs m_n \geqs  a \et cl(|s_n|) = cl(a).$$
Il existe donc un entier $q$ tel que $|s_n| < a^q$. En particulier, on devra avoir  $a  \leqs m_n < a^q$ de sorte que l'on a $q\geqs 2$. On choisit une suite pour laquelle l'entier $q$ est le plus petit possible. On aura ainsi
$$a^{q-1} \leqs |s_n| < a^q.$$
On sait que l'on peut extraire une sous-suite  pour laquelle $p_n$ est strictement d\Žcroissante ou bien $r_n$ est strictement d\Žcroissante. Les deux cas se traitent de la m\me mani\re. Supposons que ce soit $p_n$, alors deux cas peuvent se pr\Žsenter : ou bien  $cl(p_n) = cl(|s_n|)$ et $p_n < a^{q-1}$, ou bien  $cl(p_n) < cl(|s_n|)$; dans les deux cas, il y a contradiction. \qed
 
\

\noi {\bf Cela prouve que l'ensemble $\mathbf C$ des produits est bien ordonn\Ž !}

\

\su{Lemme B} {\sl Pour chaque $t\in C$, il n'y a qu'un nombre fini de s\Žquences $s\in B$ telles que $|s| = t$.} Autrement dit, chaque produit n'a qu'un nombre fini de repr\Žsentants.

\su{D\Žmonstration} On raisonne par l'absurde. On suppose qu'il y a des produits qui poss\dent une infinit\Ž de repr\Žsentants. L'ensemble des produits, $C$, \Žtant bien ordonn\Ž, soit $t$ le plus petit de ces produits qui ont une infnit\Ž de repr\Žsentants, $s_1, s_2, \dots$, o\ $s_n = (u_n, v_n, \dots, w_n)$. Il vient $|s_n| = u_n|(v_n,\dots,w_n)|$. On extrait de la suite des $u_n$, une sous-suite croissante, $u_{\sigma(1)} \leqs u_{\sigma(2)} \leqs u_{\sigma(3)}\leqs \dots$. La suite correspondante des $z_n = |(v_n,\dots,w_n)|$ est  d\Žcroissante et devient donc stationnaire, \ˆ partir d'un certain rang, $z_n = z_{n+1} = z_{n+2} = \cdots$. Ainsi, le produit $z_n$, strictement plus petit que $t$, poss\de une infinit\Ž de repr\Žsentants. Une contradiction. \qed

\

\noi {\bf  Cela ach\ve la d\Žmonstration du lemme de H. B. Neumann.}

\

\head{Un cas particulier} 

\

Soient $(G, \leqs)$ un groupe commutatif totalement ordonn\Ž et $A$ une partie de $G$ form\Že d'\Žl\Žments strictement positifs, bien ordonn\Že pour la relation $\geqs$. L'ensemble $C$ de toutes les sommes d'un nombre fini d'\Žl\Žments de $A$ est lui-m\me bien ordonn\Ž pour la relation $\geqs$. De plus, pour chaque  somme $t\in C$, il n'y a qu'un nombre fini de repr\Žsentants de $t$ dans l'ensemble $B$ des s\Žquences.

\

\noi Ce r\Žsultat interviendra dans un expos\Ž sur les nombres de Cuesta-Conway, \ˆ venir.

\

\centerline{ $*$ \ $*$ \ $*$}

\

{\bf Je tiens \ˆ remercier, bien vivement, mon coll\gue et ami, Charles Helou, pour m'avoir procur\Ž une copie de l'article de Harkelroad et Gonshor.}

\

 \centerline{$*$  \ $*$ \ $*$}
 
 \
 
\head{Une version en anglais. Translation into English}

\

\head{\bf A COMBINATORIAL  LEMMA OF H. B. NEUMANN}

\
\

\hfill{\it And sometimes add, but oft'ner take away.} 

\hfill{Dryden, after Boileau}

\

\rm

\

\

\su{Abstract} We give a notably simpler and shorter proof of H. B. Neumann's result which is stated, cursorly, like this. For any well-ordered subset, $A$, of a totally ordered semigroup, the set of products of any finite number of elements of $A$ is itself well-ordered. Moreover, for each $t$, there are only a finite number of such products equal to $t$.

\

We find in Neumann, [ 3, p.209 ], a result that Gonshor took over, in [ 1, p.52 ], under the title {\sl Combinatorial lemma on semigroups}, and which he states as follows.

\

\noi Let $A$ be a well-ordered subset of the set of strictly positive elements of a totally ordered semigroup. Then the set of all sums, $a_1 + a_2 +\dots + a_q$, of any finite number of elements of $A$ is well-ordered and, moreover, for each given $t$, there are only a finite number of such sums which are equal to $t$.

\

\noi  Gonshor's proof is not easy to follow, nor is Neumann's. They both abound with indices, sub-indices and various accents! Gonshor implies, without explicitly saying so in his statement, that the semigroup is commutative! Indeed, the operation is denoted there, +, and he uses the commutativity in case 4 of his proof.

\

\noi While Neumann establishes the result for any totally ordered semigroup, not necessarily commutative! Its operation is multiplication. His proof occupies several pages, (p.204 -209).

\

\lopar  Reading the paper [ 2 ] by Harkelroad and Gonshor, one will find a probable confirmation of the fact that Gonshor intended to use a commutative semigroup, alias abelian.\}

\

\noi What complicates matters a little bit is that {\sc Neumann} and {\sc Gonshor} do not mean quite the same thing when speaking of a totally ordered semigroup.

\

\noi We will stick to Neumann's choice, the most general one, and give a simpler proof of his lemma. But first, these few preliminaries in order to fix notations and vocabulary and to clarify the meaning of the statements!

\

\head{Preperation}

\

\head{What is a totally ordered semigroup according to Neumann}

\

\noi By definition, a totally ordered semigroup is endowed with a binary operation, the product, and a total order relation, $\leqs$.
The product is associative, in other words, we have

$a(bc) = (ab)c$.

\noi The product and the total order are compatible. In other words, when we have $a < b$, we will have $ac < bc$ and $ca < cb$, for all $c$. 

\noi The smallest of all elements is a unit element, 1.

\

We easily deduce the following properties of right and left cancellation.

\noi If $ac = bc$ then $a = b$.

\noi Similarly, if $ca = cb$ then $a = b$.

\noi So $a=b \iff ac=bc \iff ca=cb$.

\noi Thus, $ac<bc \iff ca<cb \iff a<b$, for all $a, b, c$. 

\noi Finally, $a \leqs b \iff ac \leqs bc \iff ca \leqs cb$.

\noi In particular, $a > 1 \iff aa > a$,

\noi and $aa = 1 \iff a = 1$.

\

\lopar By having added a unit element, 1, smaller than all the others, we considerably simplify {\sc Neumann}'s considerations.\}

\

\

\su{Archimedean classes} We introduce the following binary relation between the elements of the semigroup:

$a\sim b \iff \ \text{there exists}\ m, n\in \N, \ \text{such that}\ b \leqs a^m \ \text{and} \ a \leqs b^n$.

\noi We verify that it is an equivalence relation whose classes are intervals which form a partition of the semigroup. For each element $a$ of the semigroup, we designate by $cl(a)$ the (archimedean) equivalence class to which it belongs, in other words,

$cl(a) = \{b : b \sim a\}$,

$cl(1) = \{1\}$, in particular.

\noi The set of Archimedean classes is totally ordered by the following relation

$cl(a)<cl(b) \iff a^n <b, \ \text{for all} \ n\in \N$.

\

\su{Remark} The class of a product, $p = ab\dots c\dots d$, is equal to the class of the largest factor! Indeed, suppose that all the factors are $\leqs  c$. Then $c \leqs p$  because all the factors are $\geqs 1$, on the one hand. On the other hand $p \leqs  c^n$  where $n$ is the number of factors.

\

\head{Well-ordered sets. Reminders}

\

By definition, a well-ordered set is an ordered set in which any non-empty subset has an element smaller than all the others! An ordered set, $E$, is well-ordered if and only if there is no strictly decreasing sequence in $E$. This is a well-known characterization of well-ordered sets. Here is another, perhaps less often cited.

\

Let $E$ be an ordered set. The following two statements are equivalent.

1. The set $E$ is well ordered.

2. From each sequence, $x_1, x_2, \dots$, in $E$ we can extract an increasing subsequence, $x_{\sigma(1)} \leqs x_{\sigma(2)},\dots$.

\ 

\noi {\bf Indeed}, if $E$ is well-ordered, from  the sequence $x_1, x_2, \dots$, we extract $x_{\sigma(1)}$, the smallest term of the sequence, having the smallest index, then $x_{\sigma(2)}$, the smallest term of what remains of the sequence, having the smallest index, and we start again! We obviously obtain an increasing extracted subsequence!

\

\noi Conversely, if statement 2 is satisfied, there is obviously no strictly decreasing sequence in $E$.\qed

\

\head{Strings}

\

\noi Let $A$ be a well-ordered subset of a totally ordered semigroup such that $1\notin A$. We will call {\bf string}  any $q$-tuple, $s = (a, b, \dots , c)$, where $a,b,\dots, c$, is a finite sequence of $q$ elements of $A$. We say that the string is of length $q$, and we set $|(a,b,\dots,c)| = ab\dots c$, the product.

\

\noi Beware, the semigroup not necessarily being commutative, we cannot interchange the terms of a product! Of course, two distinct strings, of equal or different lengths, can have the same product!

\

\noi It is important to distinguish a sequence, $s$, from its product, $|s|$. We denote by $B$ the set of all strings and by $C = \{|s| : s \in B\}$ the set of all products, so that we have
$$A\inc B\ , \ A\inc  C \ , \ \text{and an map}\ , B \to C \ ,  \ s\mapsto |s|.$$
We say that the string, $s$,  is a {\bf representative} of the product $|s|$.

\

We will establish the following result, by contradiction. The first part of {\sc Neumann}'s combinatorial lemma will follow. We will try to avoid excessive use of indices!

 \
 
 \su{Lemma A} {\sl There is no sequence of strings}, 
 $$s_1,s_2,\dots, s_n,\dots, \ \text{such that}$$ 
 $$ |s_1| > |s_2| > \dots  > |s_n| > \cdots$$
\su{Proof}. Assume that there are such sequences. Let $\cal S$ be the set of all these sequences. We take one of them. For 
$s_n = (a_n,b_n,\dots,c_n,\dots, d_n)$, we set
$$m_n = \max\{a_n,b_n,\dots,c_n,\dots, d_n\}.$$
We get
$$|s_n| = a_nb_n\dots c_n\dots d_n \geqs m_n \ \text{and}\ cl(|s_n|) = cl(m_n),$$
according to the remark made above. The $m_n$'s are elements of the well-ordered set, $A$. Let $m$ be the smallest of these $m_n$'s. The class $cl(m)$ is the smallest of the $cl(|s_n|)$ classes. We have $cl(|s_n|) \geqs cl(|s_{n+1}|)$ since $|s_n| > |s_{n+1}|$. We will therefore have, from a certain point on,
$$cl(m) = cl(|s_n|) = cl(|s_{n+1}|) = cl(|s_{n+2}|) = \cdots,$$
in other words, all products $|s_n|$ will be in the same Archimedean class, $cl(m)$. By deleting some of the first terms of the sequence, we are led to the case where we have
$$cl(m) = cl(|s_1|) = cl(|s_2|) = cl(|s_3|) = \cdots .$$
Each  $|s_n|$ is of the form:

$|s_n| = p_nm_nr_n$ where $p_n$ and $r_n$ are themselves products or else,
exceptionally, $p_n = 1$ or $r_n = 1$. We have $m_n\geqs m$ for all $n$. We can therefore extract a strictly decreasing subsequence of $p_n$, or else a strictly decreasing subsequence of $r_n$ because, otherwise, we would have $p_{n+1} \geqs p_n$ and $r_{n+1} \geqs r_n$, from a certain point on,  therefore $|s_{n+1}| \geqs |s_n|$.

\

\noi Among all the sequences belonging to the set $\cal S$, we take one for which the class $cl(m)$ is the smallest possible! Then we denote by $a$ the smallest of the elements of $A$ which belong to $cl(m)$. We get
$$|s_n| \geqs m_n\geqs a \ \text{and} \ cl(|s_n|) = cl(a).$$
So there exists an integer $q$ such that $|s_n| < a^q$. In particular, we should have $a \leqs m_n <a^q$ so that $q\geqs 2$. We choose a sequence for which the integer $q$ is the smallest possible. We will thus have $$a^{q-1}\leqs |s_n| < a^q.$$
We know that we can extract a subsequence for which $p_n$ is strictly decreasing or $r_n$ is strictly decreasing. Both cases are treated in the same way. Suppose it is $p_n$, then either $cl(p_n) = cl(|s_n|)$ and $p_n < a^{q-1}$, or else $cl(p_n) < cl(|s_n|)$; in both cases there is a contradiction.\qed
 
 \
 
\noi  {\bf This proves that the set $C$ of products is well ordered!}

\

\su{Lemma B} {\sl For each $t \in C$, there is only a finite number of strings, $s \in B$ such that $|s| = t$.} In other words, each product has only a finite number of representatives.

\

\su{Proof} The proof is by contradiction. Suppose there are products that have an infinite number of representatives. The set of all products, $C$, being well-ordered, let $t$ be the smallest of these products which have an infinity of representatives, $s_1,s_2,\dots$, where $s_n = (u_n,v_n,\dots, w_n)$. We get  $|s_n| = u_n|(v_n, \dots, w_n)|$. From the sequence of $u_n$'s, we extract an increasing subsequence, $u_{\sigma(1)} \leqs u_{\sigma(2)} \leqs u_{\sigma(3)}\leqs \dots$.  The corresponding sequence of $z_n = |(v_n, \dots, w_n)|$ is decreasing and therefore becomes stationary, from a certain point on, $z_n = z_{n+1} = z_{n+2} = \cdots .$ Thus, the product $z_n$, strictly smaller than $t$, has an infinity of representatives. A contradiction.\qed

\

\noi {\bf This completes the proof of  H. B. Neumann's lemma.}
 
 \
 
 \head{A special case}
 
 \
 
Let $(G, \leqs$) be a totally ordered commutative group and $A$ be a subset of $G$ well-ordered for the relation $\geqs$, formed of strictly positive elements. The set $C$ of all the sums of a finite number of elements of $A$ is itself well-ordered for the relation $\geqs$. Moreover, for each  $t \in
 C$, there is only a finite number of representatives of $t$ in the set $B$ of strings.
 
 \
 
This result will appear in a paper about Cuesta-Conway numbers, to come.

 \
 
 \centerline{$*$ \ $*$ \ $*$}
 
 \
 
{\bf With my warmest thanks to my colleague and friend, Charles Helou, for providing me with a copy of the paper by Harkelroad and Gonshor.}
 
\

\centerline{$*$ \ $*$ \ $*$}

\

\head{R\Žf\Žrences}

\

\noi[ 1 ] {\sc H. GONSHOR}, {\sl An introduction to the theory of surreal numbers}, London Mathematical Society Lecture Note Series. 110, Cambridge University Press, 1987, Re-issued in digitally printed version 2008, ii + 192 pages.

\

\noi[ 2 ] {\sc L. HARKELROAD; H. GONSHOR},
{\sl The ordinality of additively generated sets},
Algebra Univers., {\bf 27}, (4) , (1990) 507-510.

\

\noi[ 3 ] {\sc B. H.  NEUMANN}, {\sl On ordered division rings}, Trans. Amer. Math. Soc. {\bf 66} (1949) 202-252. 

\

\

\noi Classification MSC : 20M99, 05E14, 06M99, 11B99

\

\noi Mots-clefs : semigroupe, totalement ordonn\Ž 

\

\enddocument